\documentclass{ws-jktr}
\usepackage{graphicx}
\begin{document}

\catchline{}{}{}{}{}


\title{Detecting and visualizing $3$-dimensional surgery}
\author{Stathis Antoniou}
\address{ \tiny {School of Applied Mathematical and Physical Sciences, National Technical University of Athens, Greece \newline  \textit{santoniou@math.ntua.gr}}}  
\author{Louis H.Kauffman}
\address{\tiny {Department of Mathematics, Statistics, and Computer Science, University of Illinois at Chicago, Chicago, USA \newline Department of Mechanics and Mathematics, Novosibirsk State University, Novosibirsk, Russia\newline \textit{kauffman@uic.edu}}}
\author{Sofia Lambropoulou}
\address{\tiny {School of Applied Mathematical and Physical Sciences, National Technical University of Athens, Greece \newline \textit{sofia@math.ntua.gr}}  }
\maketitle

\begin{abstract}
Topological surgery in dimension $3$ is intrinsically connected with the classification of $3$-manifolds and with patterns of natural phenomena. In this expository paper, we present two different approaches for understanding and visualizing the process of $3$-dimensional surgery.
In the first approach, we view the process in terms of its effect on the fundamental group. Namely, we present how $3$-dimensional surgery alters the fundamental group of the initial manifold and present ways to calculate the fundamental group of the resulting manifold. We also point out how the fundamental group can detect the topological complexity of non-trivial embeddings that produce knotting.
The second approach can only be applied for standard embeddings. For such cases, we give new visualizations for both types of $3$-dimensional surgery as different rotations of the decompactified $2$-sphere. Each rotation produces a different decomposition of the $3$-sphere which corresponds to a different visualization of the $4$-dimensional process of $3$-dimensional surgery. 
\end{abstract}

\keywords{topological surgery, framed surgery, topological process, 3-space, 3-sphere, 3-manifold, handle, decompactification, rotation, 2-sphere, stereographic projection, fundamental group, Poincar\'{e} sphere, topology change, visualization, knots, knot group, blackboard framing}
\ccode{2010 Mathematics Subject Classification: 57M05, 57M25, 57M27, 57R60, 57R65}

\section{Introduction}\label{Intro}
Topological surgery is a mathematical technique introduced by A.H.Wallace~\cite{Wal} and J.W.Milnor~\cite{Milsur} which produces new manifolds out of known ones. It has been used in the study and classification of manifolds of dimensions greater than three while also being an important topological tool in lower dimensions. Indeed, starting with $M=S^2$, $2$-dimensional surgery can produce every compact, connected, orientable $2$-manifold without boundary, see~\cite{Kosniowski,Ad}. Similarly, starting with $M=S^3$, every closed, connected, orientable $3$-manifold can be obtained by performing a finite sequence of $3$-dimensional surgeries, see~\cite{Wal,LickTh,Rolfsen}. Further, the surgery descriptions of two homeomorphic $3$-manifolds are related via the Kirby calculus, see~\cite{Kirby,Rolfsen}. 

But apart from being a useful mathematical tool, topological surgery in dimensions $1$ and $2$ is also very well suited for describing the topological changes occurring in many natural processes, see~\cite{SS1,SS2}. The described dynamics have also been exhibited by the trajectories of Lotka-{}-Volterra type dynamical systems, which are used in the physical and biological sciences, see~\cite{SS5,SS6,N7}. And, more recently, $3$-dimensional surgery has been proposed for describing cosmic phenomena such as wormholes and cosmic string black holes, see~\cite{BHsurg, SS3,SS4}.

Topological surgery uses the fact that the following two $m$-manifolds have the same boundary: $\partial (S^{n} \times D^{m-n})= \partial (D^{n+1} \times S^{m-n-1})=S^{n} \times S^{m-n-1}$. Its process removes an embedding of $S^{n} \times D^{m-n}$ (an $(m-n)$-thickening of the $n$-sphere) and glues back $D^{n+1} \times S^{m-n-1}$ (an $(n+1)$-thickening of ${(m-n-1)}$-sphere) along the common boundary, see Definition~\ref{surgery}.

In dimensions $1$ and $2$ this process can be easily understood and visualized in $3$-space as it describes removing and gluing back segments or surfaces. However, in dimension $3$ the process becomes more complex, as the additional dimension leaves room for the appearance of knotting when non-trivial embeddings are used. Moreover, the process requires four dimensions in order to be visualized. 

In this paper we present how to detect the complexity of surgery via the fundamental group and, for the case of standard embeddings, we propose a new visualization of $3$-dimensional surgery. The first approach is understanding $3$-dimensional surgery by determining the fundamental group of the resulting manifold. This approach is presented for both types of $3$-dimensional surgery, namely $3$-dimensional $0$- and $3$-dimensional $1$-surgery. For the case of $3$-dimensional $1$-surgery the presence of possible knotting during the process makes its visualization harder,
since the resulting $3$-manifolds are involved with the complexity of the knot complement. However, this complexity can be detected by the fundamental group, as we can describe the framing longitude in terms of the generators of the fundamental group of the knot complement  $S^3 \setminus N(K)$ and understand the process of $3$-dimensional $1$-surgery as the process which collapses this longitude in the fundamental group of the resulting manifold. 
On the other hand, when the standard embedding is used, we can produce a simple visualization of this $4$-dimensional process. Namely, the second approach presents a way to visualize the elementary steps of both types of $3$-dimensional surgery as rotations of the plane. As we will see, each rotation produces a different decomposition of the decompactified $3$-sphere which corresponds to a visualization of the local process of each type of elementary $3$-dimensional surgery. 

It also worth adding that, except from the two approaches presented here, there are  other ways of understanding the non-trivial embeddings of $3$-dimensional $1$-surgery. For example, the whole process can be seen as happening within a $4$-dimensional handle. This approach is discussed in~\cite{SS4}. Further details on this   perspective of surgery on framed knot and Kirby calculus can be found in~\cite{StiGompf}.

The paper is divided in four main parts: in Section~\ref{MDProcessBig} we define the processes of topological surgery. Then, in Section~\ref{EmbTopoChange3D} we discuss the topology change induced by a sequence of $3$-dimensional surgeries on a $3$-manifold $M$ and point out the role of the embedding in the case of $3$-dimensional $1$-surgery. In Section~\ref{Fundappendix} we present our first approach, which uses the fundamental group of the resulting manifold for understanding $3$-dimensional surgery. In order to analyze the case of knotted surgery curves in $3$-dimensional $1$-surgery, we also present the blackboard framing of a knot as well as the knot group. Our second approach is discussed in Section~\ref{TrivialLongi}, where we point out how the stereographic projection of the $m$-sphere can be used in order to visualize $m$-dimensional surgery in one dimension lower and we present the visualization of the elementary steps of each type of $3$-dimensional surgery via rotations of the stereographic projection of the $2$-sphere.

\section{The process of topological surgery} \label{MDProcessBig}
In Section~\ref{MDProcess} we define the process of $m$-dimensional $n$-surgery while in Sections~\ref{Types2D} and~\ref{Types3D} we present the types of $2$ and $3$-dimensional surgery.

\subsection{The process of $m$-dimensional $n$-surgery} \label{MDProcess}
\begin{definition} \label{surgery} \rm An \textbf{$m$-dimensional $n$-surgery} is the topological process of creating a new $m$-manifold $M'$ out of a given $m$-manifold $M$ by removing a framed $n$-embedding $h:S^n\times D^{m-n}\hookrightarrow  M$, and replacing it with $D^{n+1}\times S^{m-n-1}$, using the `gluing' homeomorphism $h$ along the common boundary $S^n\times S^{m-n-1}$. Namely, and denoting surgery by  $\chi$:
$$M' = \chi(M) = \overline{M\setminus h(S^n\times D^{m-n})} \cup_{h|_{S^n\times S^{m-n-1}}} (D^{n+1}\times S^{m-n-1}). $$
The resulting manifold $M'$ may or may not be homeomorphic to $M$. Note that from the definition, we must have $n+1 \leq m$. Also, the horizontal bar in the above formula indicates the topological closure of the set underneath. For details the reader is referred to~\cite{Ra}.
\end{definition}

\subsection{Types of $2$-dimensional surgery}\label{Types2D}
In dimension $2$, the above Definition gives two types of surgery. For $m=2$ and $n=0$, we have the $2$-dimensional $0$-surgery, whereby two discs $S^0\times D^2$ are removed from a $2$-manifold $M$ and are replaced in the closure of the remaining manifold by a cylinder $D^1\times S^1$:
\begin{samepage} 
 \begin{center}
$\chi(M) = \overline{M\setminus h(S^0\times D^{2})} \cup_{h} (D^{1}\times S^{1})$
\end{center}  
\end{samepage} 

For $m=2$ and $n=1$, \textit{$2$-dimensional $1$-surgery} removes a cylinder $S^1\times D^1$ and glues back two discs $D^2\times S^0$. We will only consider the first type of surgery as $2$-dimensional $1$-surgery is just the reverse (dual) process of $2$-dimensional $0$-surgery.

\begin{figure}[ht!]
\begin{center}
\includegraphics[width=7cm]{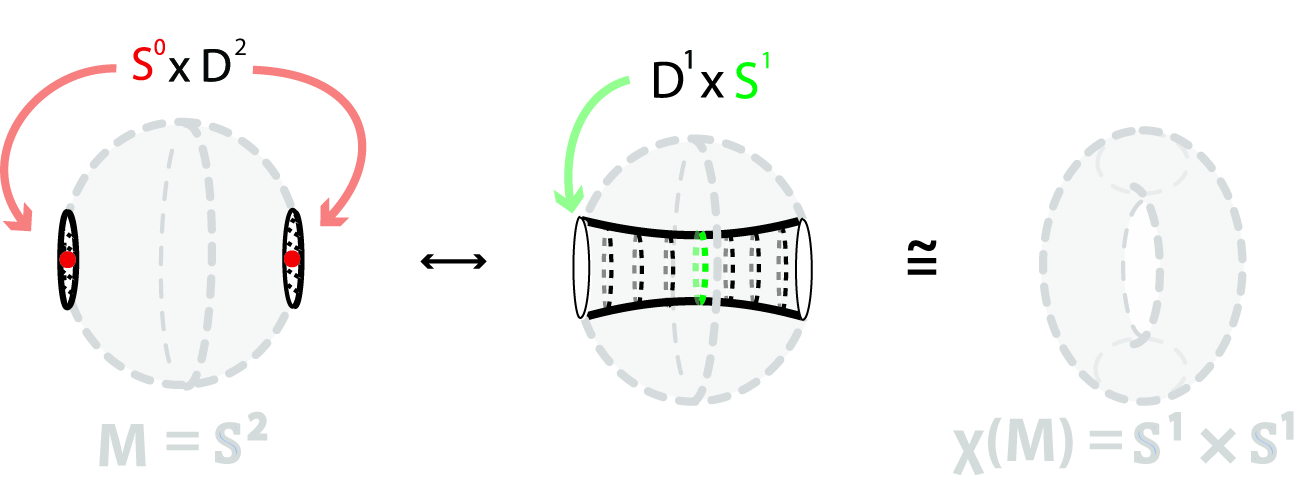}
\caption{ $2$-dimensional $0$-surgery on the sphere }
\label{2Dex}
\end{center}
\end{figure}

For example, $2$-dimensional $0$-surgery on the $M=S^2$ produces the torus $\chi(M)=S^1 \times S^1$, see Fig.~\ref{2Dex}.

\subsection{Types of $3$-dimensional surgery}\label{Types3D}
Moving up to dimension $3$, Definition~\ref{surgery} gives us three types of surgery. Namely, starting with a 3-manifold $M$, for $m=3$ and $n=0$, we have the \textit{ $3$-dimensional $0$-surgery}, whereby two 3-balls $S^0\times D^3$ are removed from  $M$ and are replaced in the closure of the remaining manifold by a thickened sphere $D^1\times S^2$:

\begin{samepage} 
 \begin{center}
$\chi(M) = \overline{M\setminus h(S^0\times D^{3})} \cup_{h} (D^{1}\times S^{2})$
\end{center}  
\end{samepage} 

Next, for $m=3$ and $n=2$, we have the \textit{ $3$-dimensional $2$-surgery}, which is the reverse (dual) process of $3$-dimensional $0$-surgery. 

Finally, for $m=3$ and $n=1$, we have the self-dual \textit{ $3$-dimensional $1$-surgery}, whereby a solid torus $S^1\times D^2$ is removed from $M$ and is replaced by another solid torus $D^2\times S^1$ (with the factors now reversed) via a homeomorphism $h$ of the common boundary:

\begin{samepage} 
 \begin{center}
$\chi(M) = \overline{M\setminus h(S^1\times D^{2})} \cup_{h} (D^{2}\times S^{1})$
\end{center}  
\end{samepage} 

For example, let us consider a $3$-dimensional $1$-surgery on $M=S^3$ using the standard embedding $h_{s}$. This embedding restricted to the common boundary $S^1\times S^1$ induces the standard pasting map $h_{s}$ which maps each longitude (respectively meridian) of $S^1\times D^2$ to a meridian  (respectively longitude) of $D^2\times S^1$. The operation produces the trivial lens space $L(0,1)$: $\chi(S^3) = \overline{S^3\setminus h_s(S^1\times D^2)} \cup_{h_s} (D^2 \times S^1) = (D^2 \times S^1) \cup_{h_s} (D^2 \times S^1) = S^2 \times S^1 =L(0,1)$.

\section{Topology change of $3$-dimensional surgeries}\label{EmbTopoChange3D}
Each type of $3$-dimensional surgery induces a different topology change on a $3$-manifold $M$. In Section~\ref{TopoChange3D0}, we discuss the topology change induced by a sequence of $3$-dimensional $0$-surgeries on $M$ and point out that the choice of the embedding $h$ in Definition~\ref{surgery} doesn't affect the resulting manifold. In Section~\ref{TopoChange3D1}, we discuss the topology change induced by a sequence of $3$-dimensional $1$-surgeries on $M$ where the embedding $h$ plays a crucial role and introduce the notion of `knot surgery'.

\subsection{$3$-dimensional $0$-surgery}\label{TopoChange3D0}
The result $\chi(M)$ of a $3$-dimensional $0$-surgery on a $3$-manifold $M$ is homeomorphic to the connected sum $M \# ({S^{1}\times S^{2}})$, independently of the embedding $h$, see~\cite{SS4}. 
Hence, we will consider the elementary step of $3$-dimensional $0$-surgery to be using the standard (or trivial) embedding $h_{s}$. We will use this fact in Section~\ref{TrivialLongi} where we present how to obtain the elementary step of the local process of $3$-dimensional $0$-surgery via rotation.

\subsection{$3$-dimensional $1$-surgery}\label{TopoChange3D1}
In contrast with $3$-dimensional $0$-surgery, $3$-dimensional $1$-surgery produces a much greater variety of $3$-manifolds. Indeed, as mentioned in Section~\ref{Intro}, every closed, connected, orientable $3$-manifold can be obtained by performing a finite sequence of $3$-dimensional $1$-surgeries on $M=S^3$, see~{\cite[Theorem 6]{Wal}} and~{\cite[Theorem 2]{LickTh}}. 
As previously, we will consider the elementary step of $3$-dimensional $1$-surgery to be using the standard embedding $h_{s}$. However, here, such elementary steps produce only a restricted family of $3$-manifolds. Indeed, starting from $S^3$, standard  embeddings $h_{s}$ can only produce $S^2 \times S^1$ or connected sums of $S^2 \times S^1$ while more complicated $3$-manifolds, such as the  Poincar\'{e}   homology sphere, require using a non-trivial embedding $h$. Hence, unlike $3$-dimensional $0$-surgery, the embedding $h$ plays an important role in the resulting manifold of $3$-dimensional $1$-surgery.      

As mentioned in Section~\ref{Types3D}, the standard embedding $h_{s}$ maps the  longitudes of the removed solid torus $V_1=S_1^1\times D^2$ to the meridians of solid torus $V_2=D^2\times S_2^1$ which is glued back and vice versa, hence $h_{s}(\ell_1)=m_2$ and  $h_{s}(m_1)=\ell_2$. When such embedding is used, the core and the longitude $\ell_1$ of the removed solid torus $V_1$ are both trivial loops, or unknotted circles. This fact allows us to obtain the elementary step of the local process of $3$-dimensional $1$-surgery via rotation. This visualization is presented with the visualization of $3$-dimensional $0$-surgery in Section~\ref{TrivialLongi}, where it is also shown that the visualizations of both types of $3$-dimensional surgeries are closely related as each one corresponds to a different rotation.

When using a non-trivial embedding $h$, both the core curve and the longitude of the removed solid torus $h(S^1\times D^2)$ can be knotted. Hence the process of $3$-dimensional $1$-surgery can be also described in terms of knots. We will call this process `knot surgery' in order to differentiate it from the process of $3$-dimensional $1$-surgery where  $h_{s}$ is used. Here, we can view the embedding $h(V_1)=h(S_1^1\times D^2)$ as a tubular neighbourhood $N(K)$ of knot $K$: $N(K)=K\times D^2=h(S_1^1\times D^2)$. The knot $K=h(S_1^1\times \{ 0 \})$ is the surgery curve at the core of solid torus $N(K)=h(S_1^1 \times D^2)$. On the boundary of $N(K)$, we further define the \textit{framing longitude} $\lambda \subset \partial N(K)$ with $ \lambda=h(S_1^1 \times \{ 1 \} )$, which is a parallel curve of $K$ on $\partial N(K)$, and the meridian $m_1 \subset \partial N(K)$ which bounds a disk of solid torus $N(K)$ and intersects the core $K$ transversely in a single point. 
 
A \textit{`knot surgery'} (or `framed surgery') along $K$ with framing $\lambda$ on a manifold $M$ is the process whereby $N(K)=h(V_1)$ is removed from $M$ and $V_2=D^2 \times S_2^1$ is glued along the common boundary. The interchange of factors of the `gluing' homeomorphism $h$ along $S_1^1 \times S_2^1$ can now be written as  $h(\lambda)=m_2$ and  $h(m_1)=l_2$.

The knottedness of $h$ makes the process harder to visualize. However, the manifold resulting from knot surgery can be understood by determining its fundamental group. In Section~\ref{Fundappendix}, we describe how to calculate this fundamental group by writing down a longitudinal element $\lambda$ in the fundamental group of the complement of the knot $S^3 \setminus N(K)$.

\section{Detecting 3-dimensional surgery via the fundamental group}\label{Fundappendix}
The fundamental group is one of the most significant algebraic constructions for obtaining topological information about a topological space.  It is a topological invariant: homeomorphic topological spaces have the same fundamental group. In Section~\ref{Fund3d0} we present how to determine the fundamental group of $3$-dimensional $0$-surgery. The more complicated topological changes, occurring during $3$-dimensional $1$-surgery, are analyzed in Section~\ref{Fund3d1}. 

\subsection{$3$-dimensional $0$-surgery}\label{Fund3d0}
As mentioned in Section~\ref{TopoChange3D0}, the resulting manifold of $3$-dimensional $0$-surgery is $M \# ({S^{1}\times S^{2}})$. The fundamental group of $\chi(M)$ can be characterized using the following lemma which is a consequence of the  Seifert-{}-van Kampen theorem (see for example~\cite{Munkres}):

\begin{lemma} \label{fundalem3d0} \rm  Let $m \geq 3$. Then the fundamental group of a connected sum of $m$-dimensional manifolds is the free product of the fundamental groups of the components: 
$$ \pi_1(M \# M')\cong \pi_1(M) * \pi_1(M')$$
\end{lemma}

Based on the above, a $3$-dimensional $0$-surgery on $M$ alters its fundamental group as follows: $\pi_1(\chi(M)) \cong \pi_1(M \# ({S^{1}\times S^2})) \cong \pi_1(M) * \pi_1({S^{1}\times S^2})\cong \pi_1(M) * ( \pi_1{(S^{1})} \times \pi_1{(S^{2})}) \cong \pi_1(M) * \mathbb{Z}$.

\subsection{$3$-dimensional $1$-surgery}\label{Fund3d1}
In Section~\ref{Blackboard}, we present the blackboard framing of a knot which will allow us to present the theorem determining the fundamental group of the manifold resulting from $3$-dimensional $1$-surgery. This is done in Section~\ref{FundaTr} 
where we also discuss the case of framed surgery along the unknot. Next, in Section~\ref{KG}, we describe the fundamental group of a knot and its presentation which allows to present the case where the surgery curve is knotted in Section~\ref{surgeryonK}.
 
\subsubsection{The blackboard framing}\label{Blackboard}
\begin{figure}[ht!]
\begin{center}
\includegraphics[width=12cm]{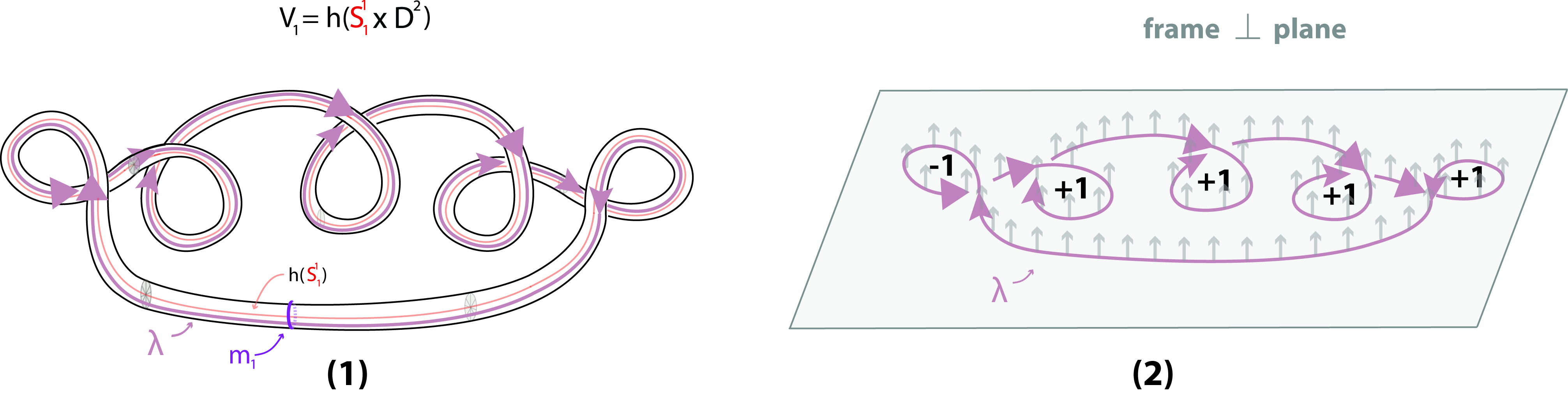}
\caption{\textbf{(1)} Isotopy of $\lambda$ \textbf{(2)} Blackboard framing of $\lambda$}
\label{3D_31_Framing}
\end{center}
\end{figure}

A framing of a knot can be also viewed as a choice of non-tangent vector at each point of the knot. The \textit{blackboard framing} of a knot is the framing where each of the vectors points in the vertical direction, perpendicular to the plane, see Fig.~\ref{3D_31_Framing}(2). The blackboard framing of a knot gives us a well-defined general rule for determining the framing of a knot diagram. Here the knot diagram is taken up to regular isotopy, namely up to Reidemeister II and III moves (see~\cite{Ad} for details on the Reidemeister moves). We use the curling in the diagram to determine the framing for an embedding corresponding to the knot, as will be explained below. Note that once we have chosen a longitude for the blackboard framing we can allow Reidemeister I moves (that might eliminate a curl) and just keep track of how the longitude now winds on the torus surface. 

\smallbreak
\begin{figure}[ht!]
\begin{center}
\includegraphics[width=12cm]{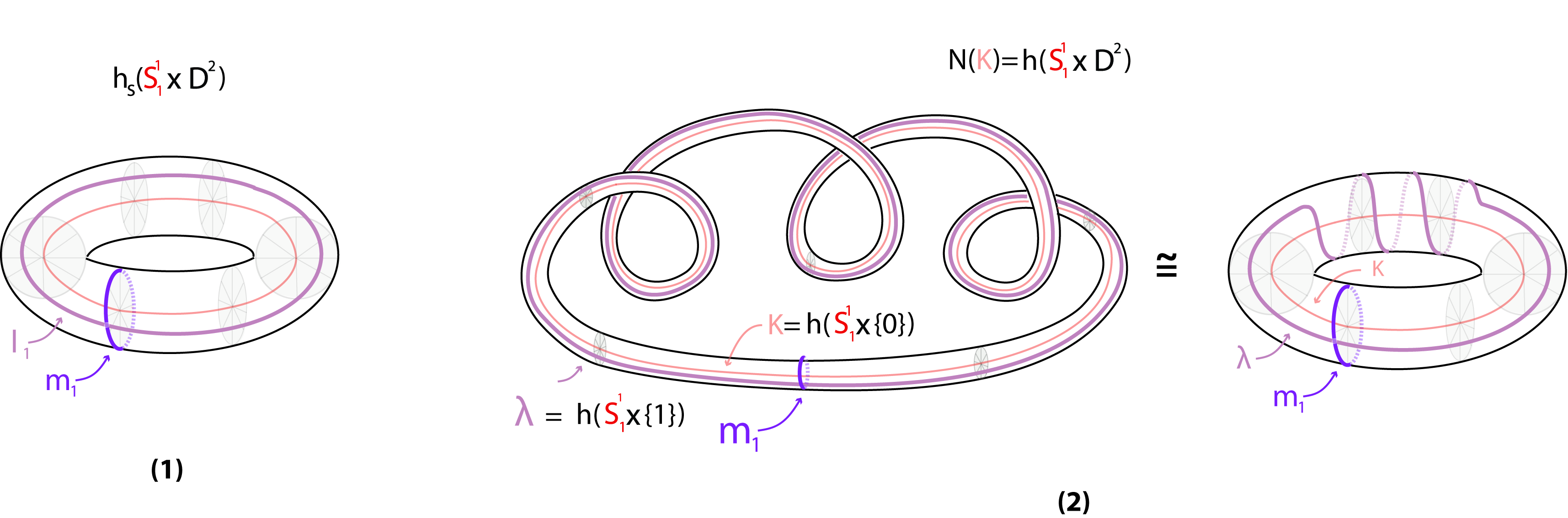}
\caption{\textbf{(1)} Longitude $l_1$ \textbf{(2)} Longitude $\lambda=l_1+3 \cdot m_1$ }
\label{3D_31_13}
\end{center}
\end{figure}

An example is shown in Fig.~\ref{3D_31_13}(2). This case  corresponds to a non-trivial embedding $N(K)=h(S_1^1 \times D^2)$ where both the knot $K$ and the longitude $\lambda$ perform three curls. As also shown in Fig.~\ref{3D_31_13}(2),  
there is an isotopic embedding of $N(K)$ where the surgery curve $K$ at the 
core of $N(K)$ is unknotted while the curls of $\lambda$ have become windings around $K$. This allows us to express $\lambda$ in terms of the unknotted longitude $l_1$ of the trivial embedding shown in Fig.~\ref{3D_31_13}(1). Namely, as $\lambda$  performs $3$ revolutions around a meridian, it can be expressed as $\lambda=l_1+3 \cdot m_1$, see Fig.~\ref{3D_31_13}(2).

More generally, if a longitude $\lambda$ performs $p$ revolutions around a meridian, it can be expressed as $\lambda=l_1+p \cdot m_1$. The induced `gluing' homeomorphism along the common boundary $S_1^1\times S_2^1$ maps each $\lambda$ of $V_1$ to a meridian of $V_2$, hence $h(l_1+p.m_1)=m_2$, while the meridians of $V_1$ are mapped to longitudes of $V_2$, hence $h(m_1)=h_{s}(m_1)=l_2$. Note that the resulting manifolds obtained by doing a $3$-dimensional $1$-surgery on $M=S^3$ using such framings on the unknot are the lens spaces $L(p,1)$. For $p=0$ we have $h(l_1)=h_s(l_1)=m_2$ and $L(0,1)=S^2 \times S^1$, which was the case presented in Section~\ref{Types3D}. For more details on lens spaces see, for example, \cite{PS}. Note that, since the multiple of the meridian $p$ is the framing number, this type of surgery is also called `framed surgery'.

\smallbreak

Recall that in Fig.~\ref{3D_31_13} the framing was $p=3$, as $\lambda$ performs $3$ revolutions. However, determining the framing of a knot diagram requires a well-defined general rule. For instance, that rule should give the same framing $p=3$ for the isotopic curve shown in Fig.~\ref{3D_31_Framing}(1). This general rule is to  take the natural framing of a knot to be its \textit{writhe}, which is the total number of positive crossings minus the total number of negative crossings. The rule for the sign of a crossing is the following: as we travel along the knot, at each crossing we consider a counterclockwise rotation of the overcrossing arc. If we reach the undercrossing arc and are pointing the same way, then the crossing is positive, see Fig.~\ref{3D_31_Framing}(2). Otherwise, the crossing is negative, see also Fig.~\ref{3D_31_Framing}(2).

Using this convention we can calculate $\lambda$ and be sure that isotopic knots will have the same framing. For instance, in Fig.~\ref{3D_31_Framing}(1), the framing number is the writhe of the knot diagram which is $p=Wr(\lambda)=4-1=3$.

\subsubsection{The fundamental group of $\chi_{\mbox{\tiny K}}(S^3)$}\label{FundaTr}
In this section, we present the theorem which characterizes the effect of knot surgery on $M=S^3$ by determining the fundamental group of the resulting manifold. We then apply it on the simple case of framed surgery along an unknotted surgery curve.

The fundamental group of the $3$-sphere $S^3$ is trivial, as any loop on it can be continuously shrunk to a point without leaving $S^3$. To examine how knot surgery alters the trivial fundamental group of $S^3$, let us consider the tubular neighborhood $N(K)$ of knot $K$. The generators of the group of $\partial N(K)$ are the longitudinal curve $\lambda$ and the meridional curve $m_1$. Note now that in $V_1=N(K)$ meridional curves bound discs while it is the specified framing longitudinal curve $\lambda$ that bounds a disc in $V_2=D^2 \times S_2^1$, since,  after surgery, the disc bounded by $m_2$ is now filling the longitude $\lambda$. Thus, $\lambda$ is made trivial in the fundamental group of $\chi_K(S^3)$. In this sense, surgery collapses $\lambda$. This statement is made precise by the following theorem which is a consequence of the Seifert-{}-van Kampen theorem (see for example~\cite{Munkres}): 
 
\begin{theorem} \label{3d1long} \rm  Let $K$ be a blackboard framed knot with longitude $\lambda \in \pi_1(S^3 \setminus N(K))$. Let  $\chi_{\mbox{\tiny K}}(S^3)$ denote the $3$-manifold obtained by surgery on $K$ with framing longitude $\lambda$.  Then we have the isomorphism:
 $$  \pi_1(\chi_{\mbox{\tiny K}}(S^3)) \cong \frac{\pi_1(S^3 \setminus N(K))}{<\lambda>} $$
where $<\lambda>$ denotes the normal subgroup generated by $\lambda \in \pi_1(S^3 \setminus N(K))$.
\end{theorem}

For a proof, the reader is referred to~\cite{Kirby,DNA}. The theorem tells us that in order to obtain the fundamental group of the resulting manifold, we have to factor out $<\lambda>$ from $\pi_1(S^3 \setminus N(K))$. 

\begin{example}  \label{ExFramedUnknot}
When the trivial embedding $h_{s}$ is used, then the `gluing' homeomorphism is $h_{s}(l_1)=m_2$, $K=the \ unknot$, $\lambda=l_1$ and $l_1$ is a trivial element in $\pi_1(S^3 \setminus N(K))$, so $<\lambda>=<0>$. In this case, we obtain the lens space $L(0,1)$ and the above formula gives us:

\begin{equation}
\begin{aligned}
\pi_1(\chi(S^3))\cong \frac{\pi_1(S^3\setminus h_s(S_1^1\times D^2))}{<\lambda>} = \frac{\pi_1(\mathring{D^2} \times S^1)}{<0>} = \frac{ \mathbb{Z}}{\{0\}} \cong\mathbb{Z}
\end{aligned} \nonumber
\end{equation}
\end{example}

\smallbreak
\begin{example}  \label{ExFramedUnknot2}
When we use a non-trivial embedding $N(K)=h(S_1^1 \times D^2)$ where the specified framing longitude $\lambda$ performs $p$ curls, the `gluing' homeomorphism is $h(\lambda)=m_2$ and, as mentioned in Section~\ref{Blackboard}, we can consider that  $K=the \ unknot$. In order to use Theorem~\ref{3d1long}, we have to find the subgroup generated by $\lambda=l_1+p \cdot m_1$ in $\pi_1(S^3 \setminus N(K))$. This subgroup is $<\lambda> = <l_1+p \cdot m_1> \cong <p \cdot m_1> \cong p \cdot <m_1>\cong p\mathbb{Z}$. In this case we obtain the lens space $L(p,1)$ and its fundamental group   is the cyclic group of order $p$:

\begin{equation}
\begin{aligned}
\pi_1(\chi(S^3))\cong \frac{\pi_1(S^3\setminus h(S_1^1\times D^2))}{<\lambda>} =\frac{\pi_1(\mathring{D^2} \times S^1)}{p\mathbb{Z}} =  { \mathbb{Z}}/{p\mathbb{Z}}  
\end{aligned} \nonumber
\end{equation}
\end{example}

As we saw in Example~\ref{ExFramedUnknot2}, if $\lambda$ is not a bounding curve in the knot complement, then we need to work out just what element $\lambda$ is in the fundamental group of the knot complement. This can be done by using one of the known presentations of the fundamental group, such as the Wirtinger presentation. A detailed presentation on the fundamental group of a knot $K$ and how we can use this presentation to determine the resulting manifold for knot surgery on $M=S^3$ along $K$ is done in next sections.

\subsubsection{The knot group}\label{KG}
The \textit{fundamental group of a knot $K$} (or the \textit{knot group}) is defined as the fundamental group of the complement of the knot in $3$-dimensional space (considered to be either $\mathbb{R}^3$ or $S^3$) with a basepoint $p$ chosen arbitrarily in the complement. The group is denoted $\pi_1(K)$ or $\pi_1(S^3 \setminus N(K))$, where $N(K)$ is a tubular neighborhood of the knot $K$. To describe this group, it is useful to have the concept of the longitude and meridian elements of the fundamental group of a knot. The longitude and the meridian are loops in the knot complement that are on the surface of a torus, the boundary of $N(K)$. 

\begin{figure}[ht!]
\begin{center}
\includegraphics[width=5.5cm]{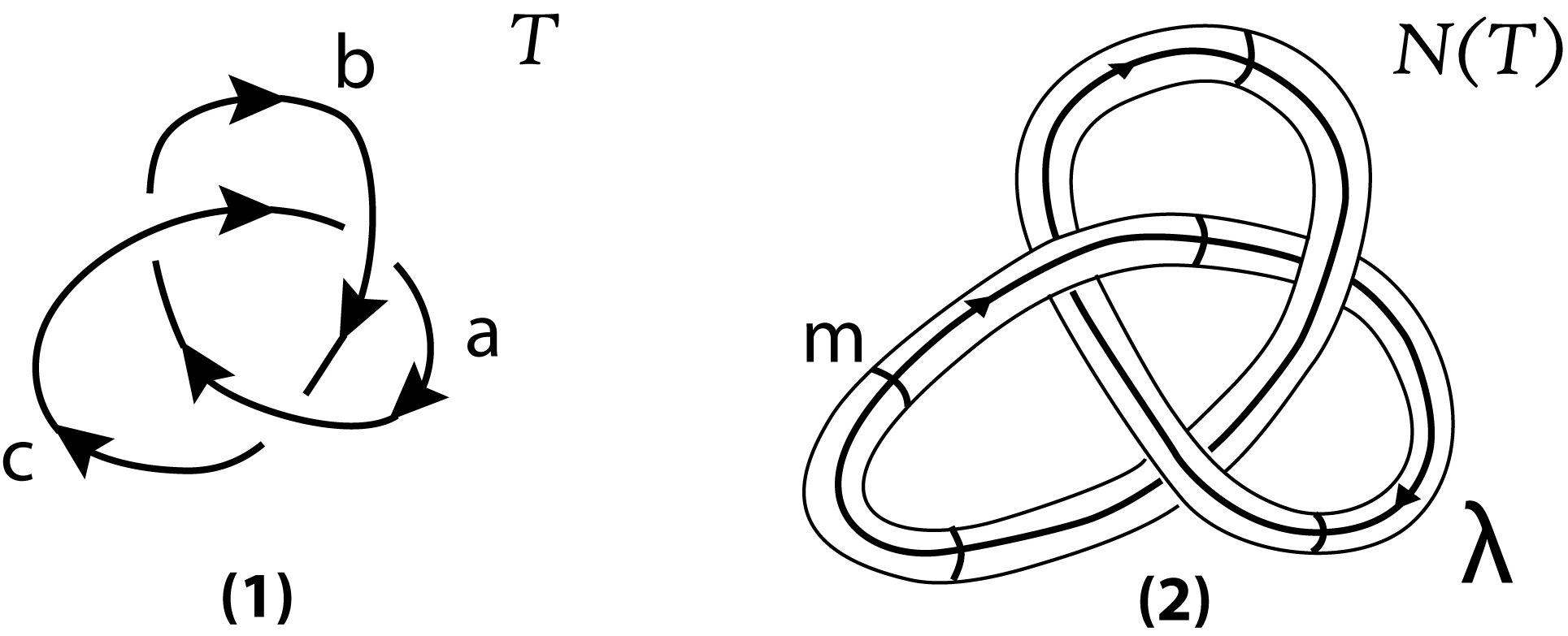}
\caption{\textbf{(1)} Trefoil knot $T$ \textbf{(2)} Tubular neighborhood $N(T)$}
\label{App_Funda_solidtrefoil}
\end{center}
\end{figure}

For the case of the trefoil knot $T$ shown in Fig.~\ref{App_Funda_solidtrefoil}(1), the meridian $m$ and the longitude $\lambda$ on the tubular neighborhood $N(T)$ are shown in Fig.~\ref{App_Funda_solidtrefoil}(2). $N(T)$ is homeomorphic to a solid torus with the knot at the core of the torus. The meridian bounds a disk in the torus, that intersects $T$ transversely in a single point. The longitude runs along the surface of the torus in parallel to $T$, and so makes a second copy of the knot out along the surface of the torus. 

\begin{figure}[ht!]
\begin{center}
\includegraphics[width=10cm]{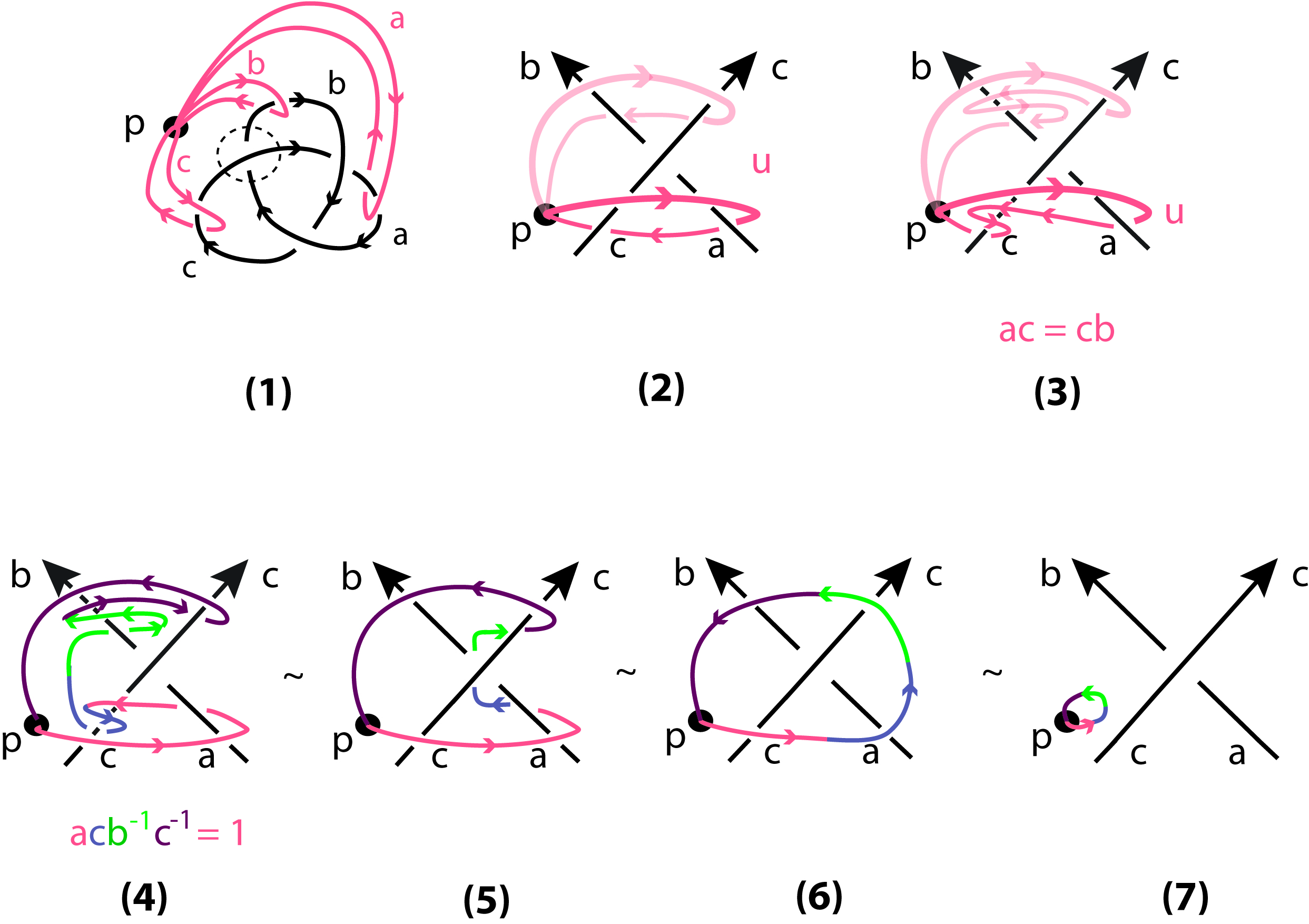}
\caption{\textbf{(1)} Generators represented as meridian loops  \textbf{(2)},\textbf{(3)} Homotopic loops \textbf{(4)},\textbf{(5)},\textbf{(6)},\textbf{(7)} Trivial curve }
\label{App_Funda_crossinglabels}
\end{center}
\end{figure}

The presentation of a knot group is generated by one meridian loop for each arc in a diagram of the knot. For the case of the trefoil, in Fig.~\ref{App_Funda_crossinglabels}(1), we illustrate the three generators $a,b,c$ (in red) which are meridian elements associated with the corresponding arcs $a,b,c$ (in black). Each crossing gives rise to a relation among those elements. For example, let  us examine the crossing of the trefoil circled in Fig.~\ref{App_Funda_crossinglabels}(1). By considering a loop $u$ in the close-up view of this crossing shown in Fig.~\ref{App_Funda_crossinglabels}(2), it is shown that $u$ wraps around arcs $a$ and $c$ but can also slide upwards to wrap around arcs $c$ and $b$. In both cases, a homotopy of loop $u$ shows that we can write $u$ as a product of the generators of the fundamental group, see Fig.~\ref{App_Funda_crossinglabels}(3). Since both homotopies describe the same loop $u$, we have $ac = cb$ which gives relation $b = c^{-1}ac$. Another way to obtain the same relation is by observing that curve $acb^{-1}c^{-1}$ contracts to a point and is therefore a trivial element of the fundamental group: $acb^{-1}c^{-1}=1$, see Fig.~\ref{App_Funda_crossinglabels}(4),(5),(6),(7).

\begin{figure}[ht!]
\begin{center}
\includegraphics[width=4cm]{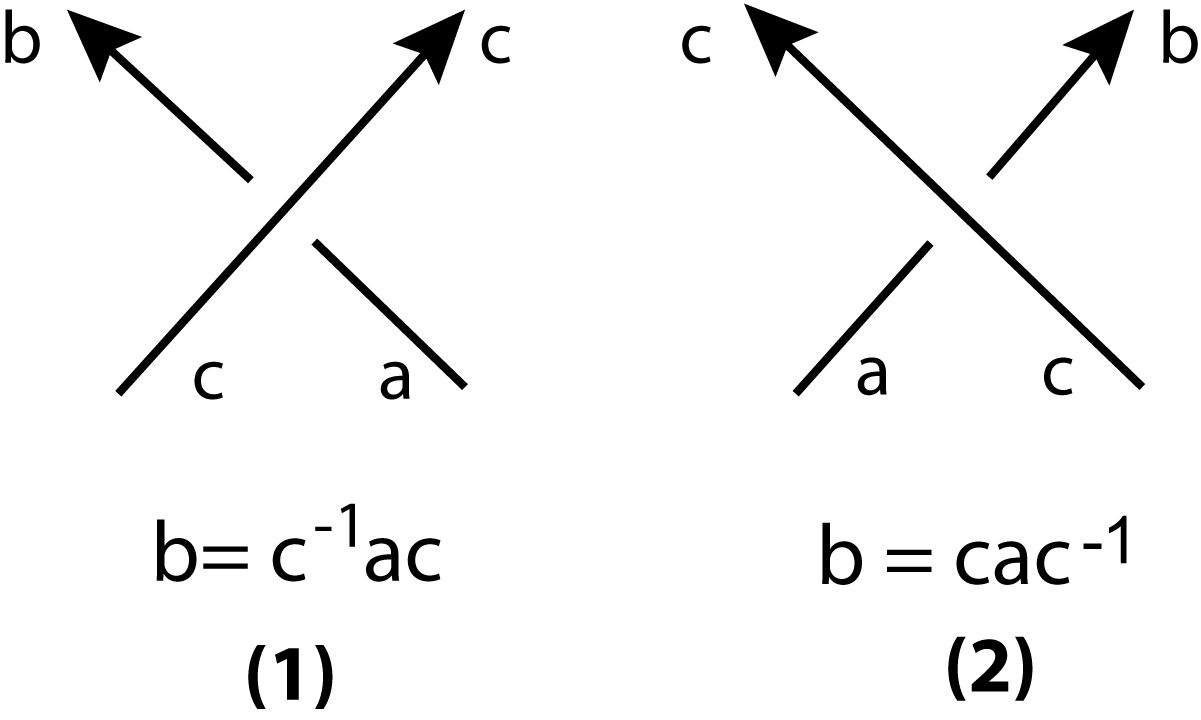}
\caption{\textbf{(1)} Positive crossing \textbf{(2)} Negative crossing}
\label{App_Funda_combi}
\end{center}
\end{figure}

Similarly, we can show that the relations obtained by the other two crossings are $a = b^{-1}cb$ and $c = a^{-1}ba$. More generally, given a diagram $D$ of an oriented knot $K$, if we label each arc of $D$, then the {\em fundamental group} of $K$ is the group whose generators are the labels of the arcs of $D$, and whose relations are the relations coming from the products of loops up to homotopy as we have just described them above. 
This presentation of the knot group is called the \textit{Wirtinger presentation} and its proof makes us of the Seifert-{}-van Kampen theorem, see for example~\cite{Rolfsen}. Hence for the trefoil knot, we have the presentation: $$ \pi_1(T)=\pi_1(S^3 \setminus N(T)) = (a,b,c \ | \ a = b^{-1}cb, b = c^{-1}ac, c = a^{-1}ba).$$ 

The fundamental group of a knot can be also defined in a combinatorial way as follows: consider a diagram of the knot and a crossing in diagram, as in Fig.~\ref{App_Funda_combi}(1) or (2), where the incoming undercrossing arc is labeled $a$, the overcrossing arc is labeled $c$ and the outgoing arc is labeled $b.$ Then write a relation in the form $b = c^{-1}ac$  for each positive crossing, as in Fig.~\ref{App_Funda_combi}(1), and a relation $b = cac^{-1}$ for each negative crossing, as in Fig.~\ref{App_Funda_combi}(2). The combinatorial approach defines the fundamental group as the group having one generator for each arc and one relation at each crossing in the diagram as we just specified them. One can  show that this group is invariant under the Reidemeister moves. This means that all diagrams of the same knot have the same fundamental group.

This combinatorial description is equivalent to the Wirtinger presentation. Indeed, see for example the relation coming from the positive crossing of Fig.~\ref{App_Funda_combi}(1) and the relation coming from homotopic loops in Fig.~\ref{App_Funda_crossinglabels}(3) or (4). However, as we will see in Section~\ref{surgeryonK}, for the purpose of doing surgery we need the topological approach,  so that we can express the longitude in terms of the generators of the fundamental group of  $S^3 \setminus N(K)$. For more details on combinatorial group theory, the reader is referred to~\cite{Sti} or~\cite{MagnusKarassSolitar}.

\subsubsection{Computing $\pi_1(\chi_{\mbox{\tiny K}}(S^3))$}\label{surgeryonK}
When the core curve $K$ of a non-trivial embedding $h(S_1^1 \times D^2)=N(K)$ is knotted, one cannot express $\lambda$ in terms of trivial longitudes and meridians, as was the case in Examples~\ref{ExFramedUnknot} and~\ref{ExFramedUnknot2}. In general, in order to compute the fundamental group of a $3$-manifold that is obtained by doing surgery on a blackboard framed knot $K$, we have to describe first how to write down a longitudinal element $\lambda$ in the fundamental group of the knot complement $S^3 \setminus N(K)$. 

\begin{figure}[ht!]
\begin{center}
\includegraphics[width=12.5cm]{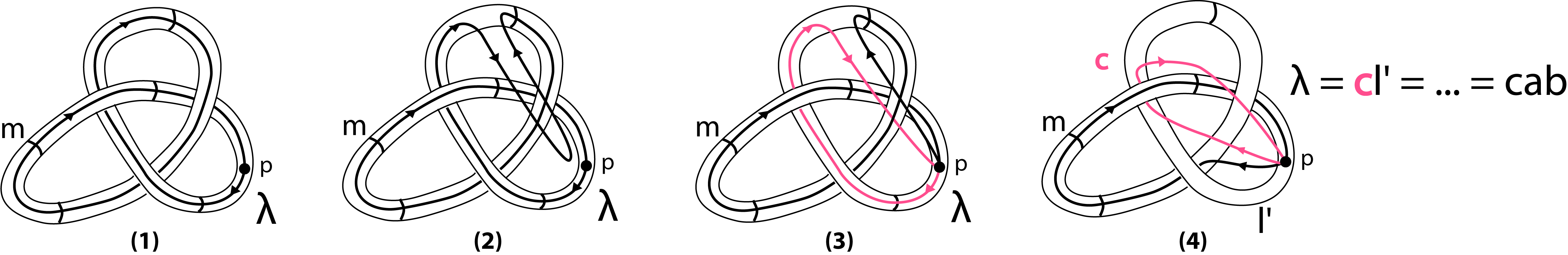}
\caption{\textbf{(1)} Longitude $\lambda$ in N(T) \textbf{(2)}\textbf{(3)}\textbf{(4)} Homotopy of $\lambda$ }
\label{App_Funda_solidtrefoil_Lambda}
\end{center}
\end{figure}

To do so, we homotope $\lambda$  to a product of the generators of $\pi_1(S^3 \setminus N(K))$ corresponding to the arcs that it underpasses. In this expression for the longitude, the elements $x$ that are passed underneath will appear either as $x$ or as $x^{-1}$   according to whether the knot is going to the right or to the left from the point of view of a traveler on the original longitude curve. Once the longitude $\lambda$ is expressed in terms of the generators of the fundamental group of  $S^3 \setminus N(K)$, we can calculate the fundamental group of $\chi_{\mbox{\tiny K}}(S^3)$ using Theorem~\ref{3d1long}.

For example, in Fig.~\ref{App_Funda_solidtrefoil_Lambda}(1) we show a trefoil knot and the longitudinal element $\lambda$ in the fundamental group running parallel alongside it. Note that, for convenience, the basepoint $p$ is on the boundary of the torus but it could be anywhere in the complement $S^3 \setminus N(K)$. Each time that $\lambda$ goes under the knot we can run a line all the way back to the base point $p$ and then back to the point where $\lambda$ comes out from underneath the knot, see Fig.~\ref{App_Funda_solidtrefoil_Lambda}(2),(3) and (4). By doing this, we have written, up to homotopy, the longitude as a product of the generators of the fundamental group that are passed under by the original longitude curve. Thus in the trefoil knot case, as shown in Fig.~\ref{App_Funda_solidtrefoil_Lambda}(4), we see that the longitude is given by $\lambda= cab$.

\begin{example} \label{TrefoilSugeries} We will now calculate the fundamental group of a $3$-manifold obtained by doing $3$-dimensional $1$-surgery on the trefoil knot  for two different projections. The first one is the simplest projection of the trefoil shown in Fig.~\ref{App_Funda_solidtrefoil_Lambda}(1). It has three positive crossings yielding a blackboard framing number of $3$. The second one has two additional negative crossings thus having a blackboard framing number of $1$, see Fig.~\ref{TrefoilP}.

\begin{figure}[ht!]
\begin{center}
\includegraphics[width=6cm]{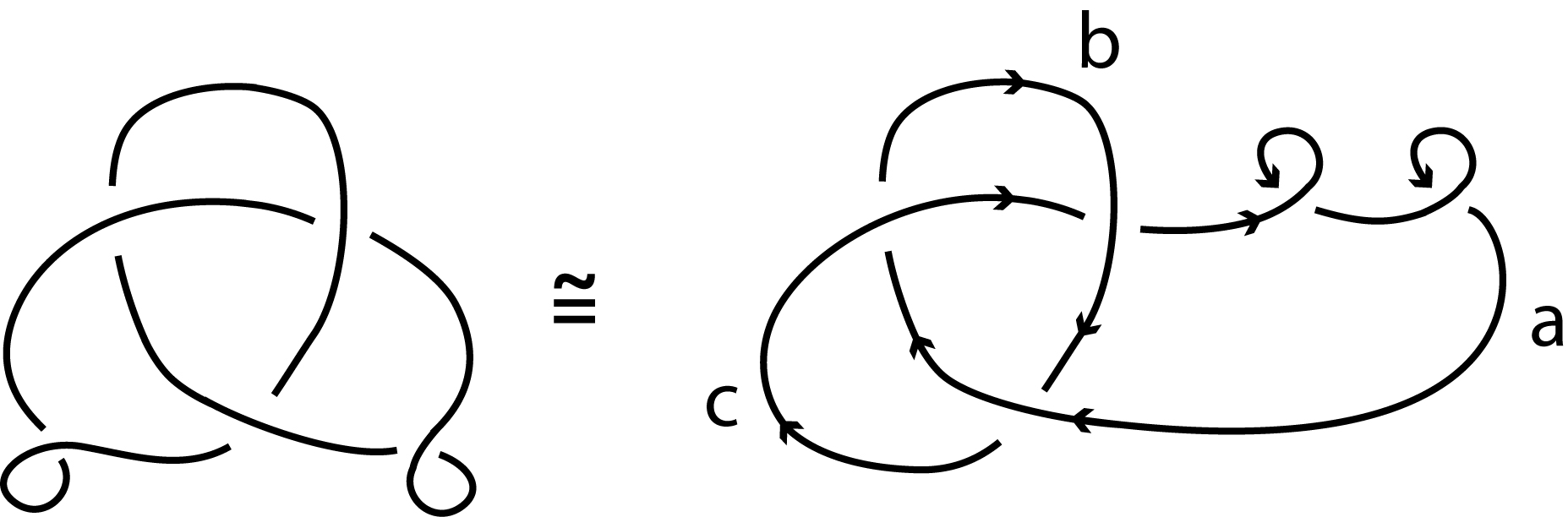}
\caption{Projection of the trefoil with total blackboard framing $1$}
\label{TrefoilP}
\end{center}
\end{figure}

As mentioned in Section~\ref{FundaTr}, surgery collapses the longitude $\lambda$, so the resulting fundamental group depends on how longitude $\lambda$ is expressed in the following relation: 

$$ \pi_1(\chi(S^3))=\frac{\pi_1(S^3\setminus N(T))}{<\lambda>} =\frac{\pi_1(T)}{<\lambda>}=(a,b,c \ | \ aba = bab, \lambda=1) $$

In the first case, by substituting $\lambda= cab$ and $c=a^{-1}ba$ to $\lambda=1$, we have $a^{-1}baab=1 \Leftrightarrow a=ba^{2}b$. Given that $aba = bab$, this implies that $a^2=baaba\Leftrightarrow a^2=babab \Leftrightarrow a^3=(ba)^3$. Notice now that  $(aba)^2=aba \cdot aba=bab \cdot bab=(ba)^3$. Thus by setting $A=a, B=ba$ and  $C=aba$ we have that  $A^3 = B^3=C^2$ and we only need to show that this is equal to $ABC$. Indeed, $ABC=a \cdot ba \cdot aba=(ba)^3$. Hence, the fundamental group of the  resulting manifold is isomorphic to the binary tetrahedral group $(A,B,C \ | \ A^3 = B^3=C^2=ABC)$ denoted $<3,3,2>$. It is also worth mentioning that the resulting manifold is isomorphic to $S^3/<3,3,2>$, the quotient of the $3$-sphere by an action of the binary tetrahedral group. For details on group actions the reader is  referred to~\cite{MilAct}. 

\begin{figure}[ht!]
\begin{center}
\includegraphics[width=5cm]{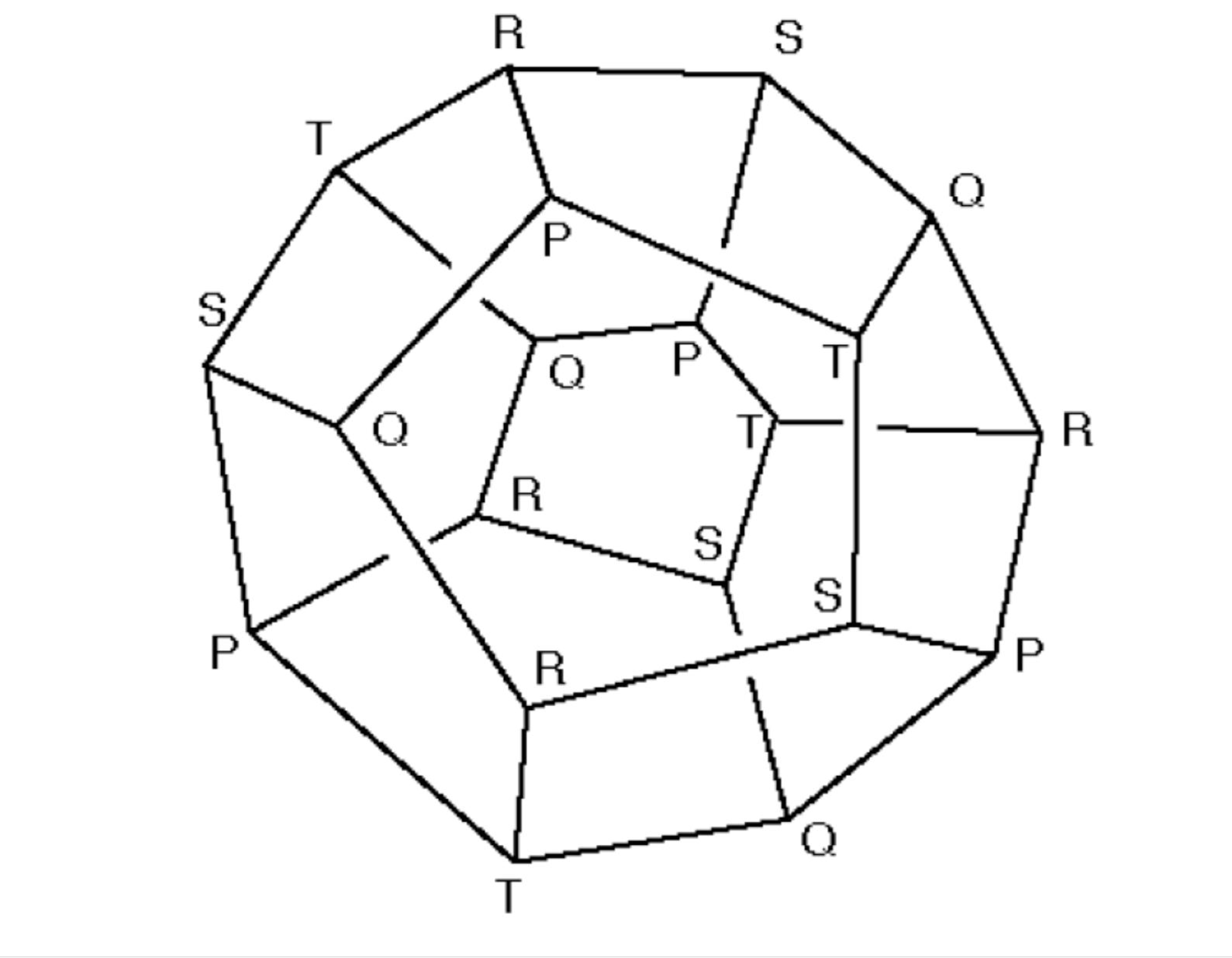}
\caption{Poincar\'{e} sphere}
\label{Poincare}
\end{center}
\end{figure}

In the second case, the longitude $\lambda$ in the projection shown in Fig.~\ref{TrefoilP} is the same as the one in Fig.~\ref{App_Funda_solidtrefoil_Lambda}(1) with two additional negative crossings along arc $a$. Hence, in this case $\lambda=caba^{-2}$. By substitution, we have $a^{-1}baaba^{-2}=1 \Leftrightarrow a^{3}=ba^{2}b$. Given that $aba = bab$, this implies that $a^3=baaba\Leftrightarrow a^4=babab \Leftrightarrow a^5=(ba)^3$.  Thus by setting $A=a, B=ba$ and  $C=aba$ we have that  $A^5 = B^3=C^2$ and   $ABC=a \cdot ba \cdot aba=(ba)^3$. The fundamental group of the  resulting manifold is isomorphic to the binary icosahedral group $(A,B,C \ | \ A^5 = B^3=C^2=ABC)$ denoted by $<5,3,2>$. The resulting manifold is isomorphic to $S^3/<5,3,2>$, the quotient of the $3$-sphere by an action of $<5,3,2>$. 

This manifold is also known as the \textit{Poincar\'{e} homology sphere}, which can be described by identifying opposite faces of a dodecahedron according to the scheme shown in Fig.~\ref{Poincare} (for more details on this identification, see~\cite{SeTh}). It can be shown from this that the Poincar\'{e} homology sphere is diffeomorphic to the link of the variety $V((2,3,5))=\{ (z_1,z_2,z_3) \ | \ z_1{^2} + z_2{^3} + z_3{^5} = 0\}$, that is, the intersection of a small $5$-sphere around $0$ with $V((2,3,5))$. From this it is not hard to see that the Poincar\'{e} homology sphere can be also obtained as a $5$-fold cyclic branched covering of $S^3$ over the trefoil knot. For more details on the different descriptions of the Poincar\'{e} homology sphere, the reader is referred to~\cite{KirScha}.

This manifold has been of great interest, and is even thought by some physicists to be the shape of the geometric universe, see~\cite{Weeks,Luminet,Levin}. See also~\cite{SS4}.
\end{example}

\section{Visualizing elementary $3$-dimensional surgery via rotation} \label{TrivialLongi}
As mentioned in Section~\ref{EmbTopoChange3D}, we consider $3$-dimensional surgeries using the standard embedding $h_{s}$ as elementary steps. In this section, we show that, in this case, both types of $3$-dimensional surgery can be visualized via rotation. To do so, we first describe how stereographic projection can be used to visualize the local process of topological surgery in one dimension lower, see Section~\ref{MDStereoBig} and then use it to visualize elementary $3$-dimensional surgeries in $\mathbb{R}^3$, see Section~\ref{3DProcess}.

\subsection{Visualizations of topological surgery using the stereographic projection} \label{MDStereoBig}
In Section~\ref{MDStereo} we present a way to visualize the initial and the final instances of $m$-dimensional surgery in $\mathbb{R}^{m}$ using stereographic projection. We then discuss the case of $m=2$ in Section~\ref{2DStereo} which will be our basic tool for the visualization of  elementary $3$-dimensional surgery  via rotation in Section~\ref{3DProcess}.

\subsubsection{Visualizing $m$-dimensional $n$-surgery in $\mathbb{R}^m$} \label{MDStereo}
Let us first mention that the two spherical thickenings involved in the process of $m$-dimensional $n$-surgery are both $m$-manifolds with boundary. Notice now, that if we glue theses two $m$-manifolds, along their common boundary using the standard  mapping $h_{s}$, we obtain the $m$-sphere which, in turn, is the boundary of the $(m+1)$-dimensional disc: $(S^n\times D^{m-n}) \cup_{h} (D^{n+1}\times S^{m-n-1}) = (\partial D^{n+1}\times D^{m-n}) \cup_{h} (D^{n+1}\times \partial D^{m-n})  = \partial (D^{n+1}\times D^{m-n})\cong \partial (D^{m+1})=S^{m}$. 

\begin{figure}[ht!]
\begin{center}
\includegraphics[width=9cm]{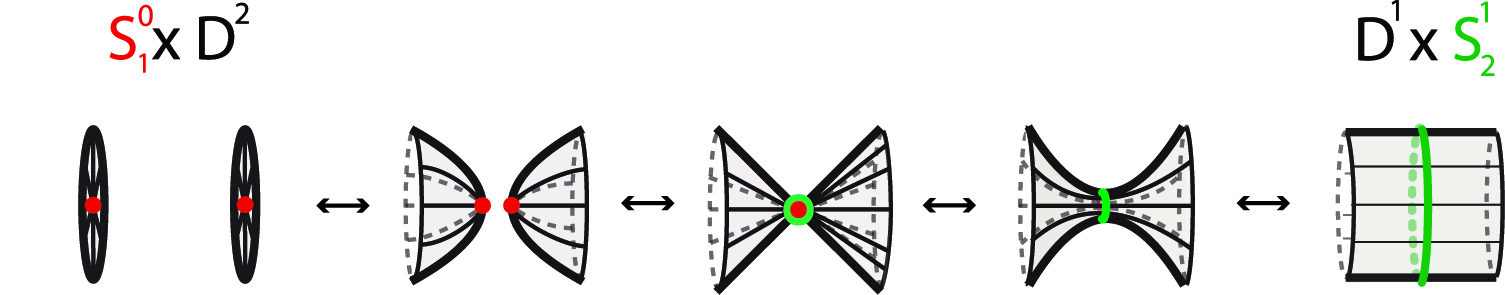}
\caption{ Continuous local process of $2$-dimensional $0$-surgery in $\mathbb{R}^{3}$ }
\label{2D_Cont}
\end{center}
\end{figure}
    
The process of surgery can be seen independently from the initial manifold $M$ as a local process which transforms $M$ into $\chi(M)$. The $(m+1)$-dimensional disc $D^{m+1} \cong D^{n+1}\times D^{m-n}$ is one dimension higher than the initial manifold $M^{m}$. This extra dimension leaves room for the process of surgery to take place continuously. The disc $D^{m+1}$ considered in its homeomorphic form $D^{n+1}\times D^{m-n}$ is an $(m+1)$-dimensional $(n+1)$-handle. The unique intersection point $D^{n+1} \cap D^{m-n}$ within $D^{n+1}\times D^{m-n}$ is called the critical point. The process of surgery is the continuous passage, within the handle $D^{n+1}\times D^{m-n}$, from boundary component $S^n\times D^{m-n} \subset  \partial (D^{n+1}\times D^{m-n})$ to its complement $D^{n+1}\times S^{m-n-1}  \subset  \partial (D^{n+1}\times D^{m-n})$ by passing through the critical point $D^{n+1} \cap D^{m-n}$. More precisely, the boundary component $S^n\times D^{m-n}$ collapse to the critical point $D^{n+1} \cap D^{m-n}$ from which the complement boundary component $D^{n+1}\times S^{m-n-1}$ emerges. For example, in dimension $2$, the two discs $S_1^0\times D^{2}$ collapse to the critical point from which the cylinder $D^{1}\times S_2^{1}$ uncollapses, see Fig.~\ref{2D_Cont}. 

Keeping in mind that gluing the two $m$-manifolds with boundary involved in the process of $m$-dimensional $n$-surgery along their common boundary gives us the $m$-sphere $S^m$, the idea of our proposed  visualization of surgery is that while $S^m$ is embedded in $\mathbb{R}^{m+1}$, it can be stereographically projected to $\mathbb{R}^{m}$. Hence, for every $m$, one can visualize the initial and the final instances of the process of $m$-surgery one dimension lower. In the following examples we deliberately did not project the intermediate instances, as this can't be done without self-intersections.

\subsubsection{Visualizing $2$-dimensional $0$-surgery in $\mathbb{R}^2$} \label{2DStereo}
\begin{figure}[ht!]
\begin{center}
\includegraphics[width=12.5cm]{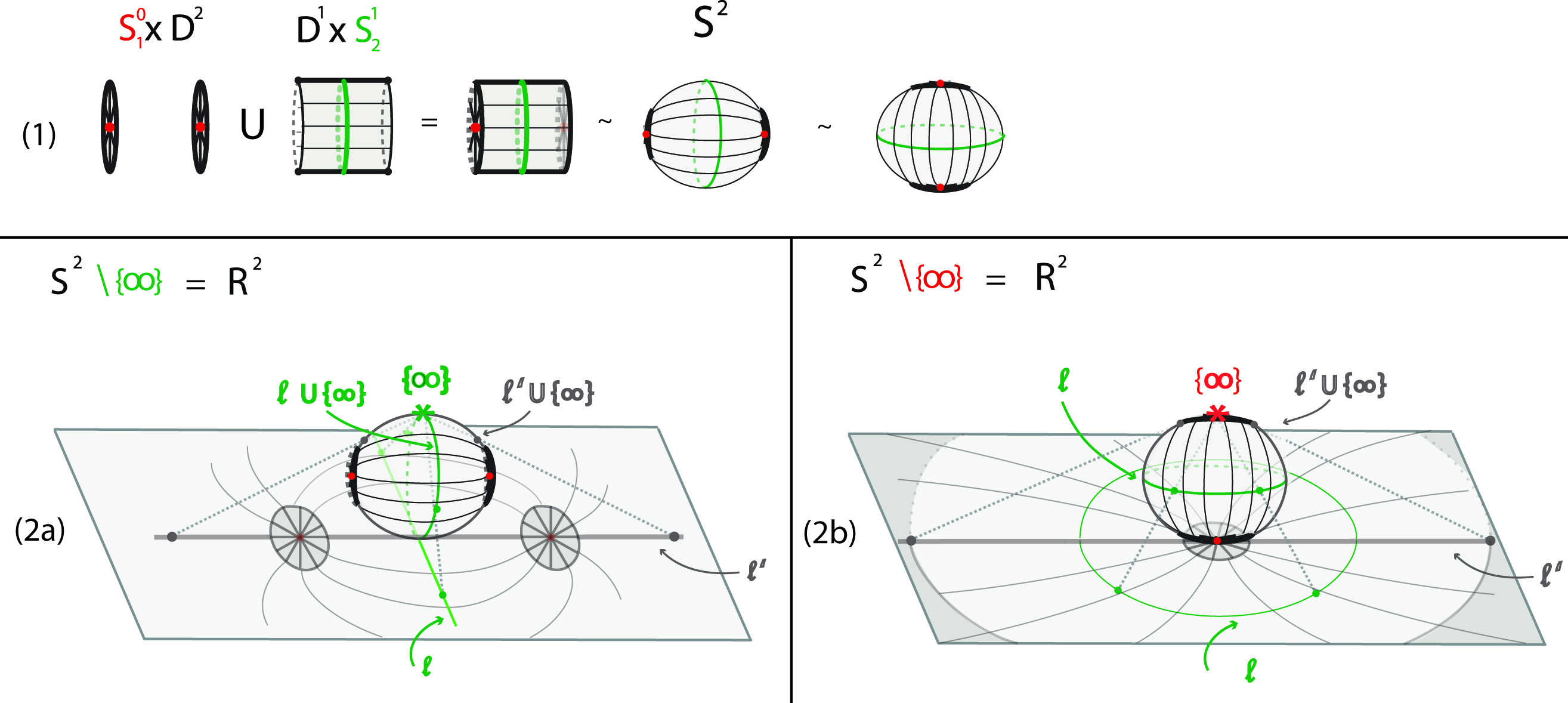}
\caption{  \textbf{(1)}  $(S_1^0\times D^{2}) \cup (D^{1}\times S_2^{1}) =S^{2}$  \newline (\textbf{{2\textsubscript{a}})} First projection of $S^2 \setminus \{\infty\}$ to $\mathbb{R}^2$ \textbf{({2\textsubscript{b}})} Second projection of $S^2 \setminus \{\infty\}$ to $\mathbb{R}^2$}
\label{2D_Proj}
\end{center}
\end{figure}

For $m=2$ and $n=0$, the initial and final instances of $2$-dimensional $0$-surgery that make up $S^2$ are shown in Fig.~\ref{2D_Proj}(1). If we remove the point at infinity, we can project the points of $S^2 \setminus \{\infty\}$ on $\mathbb{R}^2$ bijectively. We will use two different projections for two different choices for the point at infinity. The first one is shown in Fig.~\ref{2D_Proj}({2\textsubscript{a}}) where the point at infinity is a point of the core $S_2^1$ of $D^{1}\times S_2^{1}$. In this case, the two great circles $S_2^1 = \ell \cup \{\infty\}$ and $\ell' \cup \{\infty\}$ of $S^2$ are projected on the two perpendicular infinite lines $\ell$ and  $\ell'$ in ${\mathbb R}^2$. In the second one, shown in Fig.~\ref{2D_Proj}({2\textsubscript{b}}), the point at infinity is the center of one of the two discs $S_1^0\times D^{2}$. In this case the great circle $\ell' \cup \{\infty\}$ in $S^2$ is, again, projected to the  infinite line $\ell'$ in ${\mathbb R}^2$ but the great circle $S_2^1 = \ell$ is now projected to the  circle $\ell$ in ${\mathbb R}^2$. 

\begin{figure}[ht!]
\begin{center}
\includegraphics[width=12.5cm]{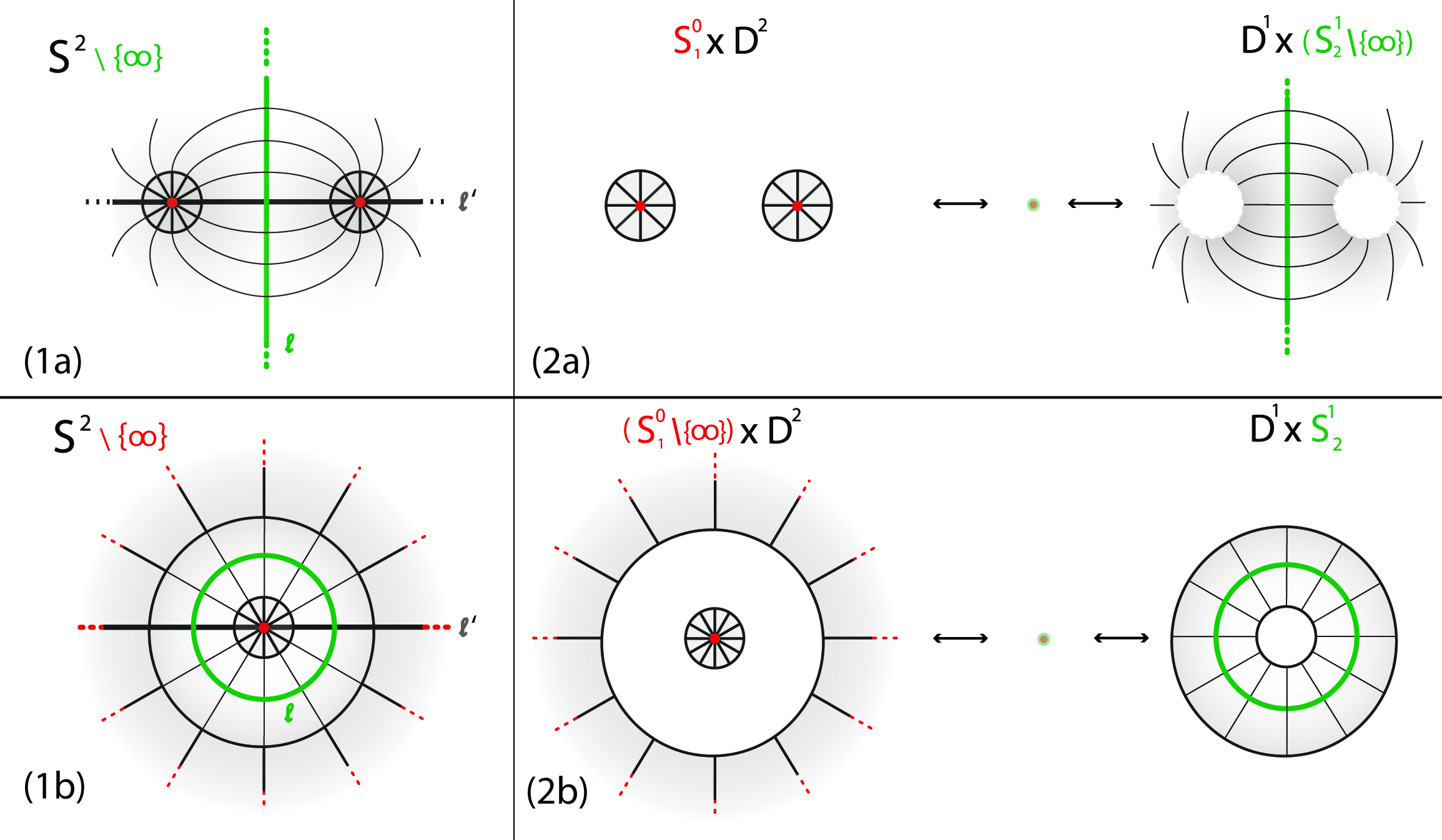}
\caption{\textbf{({1\textsubscript{a}})} First projection \textbf{({1\textsubscript{b}})} Second projection. 
\textbf{({2\textsubscript{a}})}, \textbf{({2\textsubscript{b}})} Corresponding initial and final instances of $2$-dimensional $0$-surgery in $\mathbb{R}^2$ }
\label{2D_Decomp}
\end{center}
\end{figure}

As mentioned in Section~\ref{MDStereo}, the one dimension higher of the disc $D^{m+1}$  leaves room for the process of $m$-dimensional surgery to take place continuously. For $2$-dimensional surgery, the third dimension allows the two points of the core $S_1^0$ to touch at the critical point, recall Fig.\ref{2D_Cont}. Using the two stereographic projections discussed above and shown again in Fig.~\ref{2D_Decomp}(1\textsubscript{a}) and (1\textsubscript{b}), we present in Fig.~\ref{2D_Decomp}(2\textsubscript{a}) and (2\textsubscript{b}) two different local visualizations of $2$-dimensional surgery in $\mathbb{R}^2$. Note that in Fig.~\ref{2D_Decomp}({1\textsubscript{b}}) and ({2\textsubscript{b}}), the red dashes show that all lines converge to the point at infinity which is the center of the decompactified disc and one of the points of $S_1^0$. The process of $2$-dimensional $0$-surgery starts with either one of the first instances of Fig.~\ref{2D_Decomp}(2\textsubscript{a})  and (2\textsubscript{b}). Then the centers of the two discs $S_1^0\times D^{2}$ collapse to the critical point which is shown with increased transparency to remind us that this happens in one dimension higher, see the second instances of either Fig.~\ref{2D_Decomp}(2\textsubscript{a}) or (2\textsubscript{b}). Finally the cylinder $D^{1}\times S_2^{1}$ uncollapses, as illustrated in the last instances of Fig.~\ref{2D_Decomp}(2\textsubscript{a}) and (2\textsubscript{b}). Clearly, the reverse processes provide visualizations of $2$-dimensional $1$-surgery.

\subsection{Visualizing elementary $3$-dimensional surgeries in $\mathbb{R}^3$} \label{3DProcess}


In this section we present two ways of visualizing the elementary steps of both types of $3$-dimensional surgery in $\mathbb{R}^3$ using rotations of the decompactified $2$-sphere $S^2\setminus \{\infty\}$. More precisely, in Section~\ref{Rotdecomp} we discuss the rotations of spheres and present how the rotations of $S^2\setminus \{\infty\}$ formed by the initial and  final instances of $2$-dimensional $0$-surgery in ${\mathbb R}^2$ produce the $3$-space ${\mathbb R}^3=S^3 \setminus \{\infty\}$. Then, in Section~\ref{3DProj} we detail how the rotations of the two different projections of $S^2\setminus \{\infty\}$ shown in Fig.~\ref{2D_Proj} give rise to two different decompositions of $S^3$ which correspond to two visualizations of the elementary steps of both types of $3$-dimensional surgery.

\subsubsection{Decompactification and rotations of spheres}\label{Rotdecomp}
Applying the remark of Section~\ref{MDStereo} for $m=3$ we have that the initial and final instances of all types of $3$-dimensional surgery form $S^{3}=\partial D^{4}$. Since, now, $S^3 \setminus \{\infty\}$ can be projected on ${\mathbb R}^3$ bijectively, we will present a new way of visualizing $3$-dimensional surgery in ${\mathbb R}^3$ by rotating appropriately the projections of the initial and  final instances of $2$-dimensional $0$-surgery in ${\mathbb R}^2$.

The underlying idea is that, in general, $S^n$ which is embedded in ${\mathbb R}^{n+1}$ can be obtained by a 180$^\circ$ rotation of $S^{n-1}$, which is embedded in ${\mathbb R}^n$. So, a 180$^\circ$ rotation of $S^{0}$ around an axis bisecting the interval of the two points (e.g. line $\ell$ in  Fig.~\ref{2D_Proj}({2\textsubscript{a}})) gives rise to $S^1$ (which is $\ell' \cup \{\infty\}$ in  Fig.~\ref{2D_Proj}({2\textsubscript{a}})), while a 180$^\circ$ rotation of $S^{1}$ around any diameter gives rise to $S^2$. For example, in Fig.~\ref{2D_Proj}({2\textsubscript{b}}), a 180$^\circ$ rotation of $\ell' \cup \{\infty\}$ around the north-south pole axis results in the $2$-sphere shown in the figure. Now, the creation of $S^3$ (which is embedded in ${\mathbb R}^4$) as a rotation of $S^2$ requires a fourth dimension in order to be visualized. Instead we can obtain its stereographic projection in ${\mathbb R}^3=S^3 \setminus \{\infty\}$ by rotating the stereographic projection of $S^2 \setminus \{\infty\}={\mathbb R}^2$. Indeed, a 180$^\circ$ rotation of the plane around any line in the plane gives rise to the $3$-space.

\begin{figure}[ht!]
\begin{center}
\includegraphics[width=12.5cm]{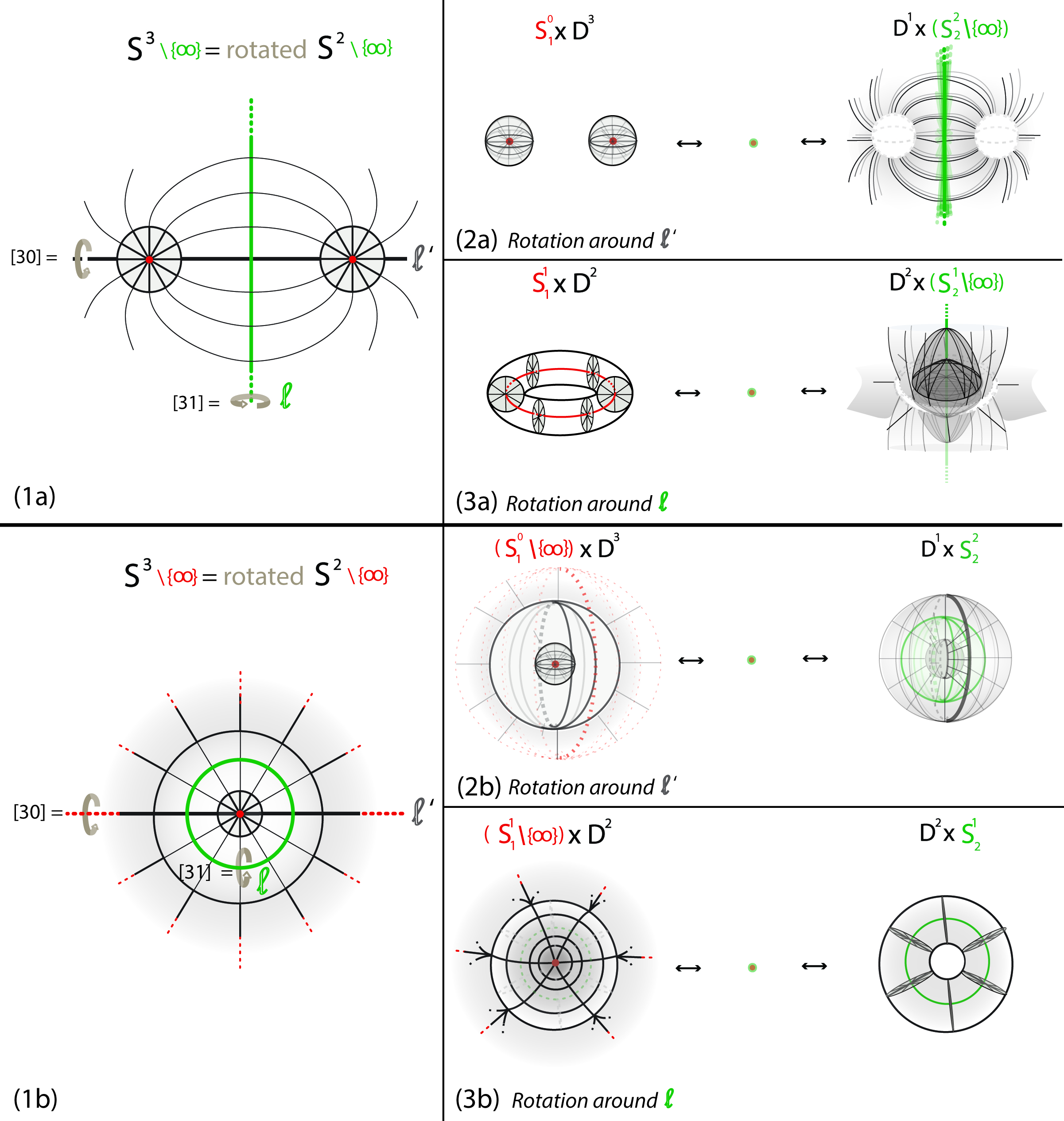}
\caption{ \textbf{({1\textsubscript{a}})},\textbf{({1\textsubscript{b}})} Representations of $S^3 \setminus \{\infty\}$ as rotations of $S^2 \setminus \{\infty\}$  using first and second projection. \textbf{({2\textsubscript{a}})},\textbf{({2\textsubscript{b}})} $3$-dimensional $0$-surgery in $\mathbb{R}^3$  using  first and second  projection. \textbf{({3\textsubscript{a}})}, \textbf{({3\textsubscript{b}})} $3$-dimensional $1$-surgery in $\mathbb{R}^3$    using first and second projection.}
\label{3D_Decomp_F}
\end{center}
\end{figure}

As we will see, each type of $3$-dimensional surgery  corresponds to a different rotation, which, in turn, corresponds to a different decomposition of $S^{3}$. As we consider here two kinds of projections of $S^2 \setminus \{\infty\}$ in ${\mathbb R}^2$, see Fig.~\ref{2D_Decomp}({1\textsubscript{a}}) and ({1\textsubscript{b}}), these give rise to two kinds of decompositions of $S^{3}$ via rotation, see Fig.~\ref{3D_Decomp_F}({1\textsubscript{a}}) and ({1\textsubscript{b}}). Each decomposition, now, leads to the visualizations of both types of $3$-dimensional surgery.
Hence, the elementary steps of both types of $3$-dimensional surgery are now visualized as rotations of the decompactified $S^2$.

\subsubsection{Projections and visualizations of elementary $3$-dimensional surgeries}\label{3DProj}
Let us start with the first projection. In Fig.~\ref{3D_Decomp_F}({1\textsubscript{a}}), we show this decompactified view in ${\mathbb R}^2$ and the two axes of rotation $\ell'$ and $\ell$. As we will see, a rotation around axis $\ell'$ induces $3$-dimensional $0$-surgery in ${\mathbb R}^3$ while a rotation around axis $\ell$ induces $3$-dimensional $1$-surgery in ${\mathbb R}^3$.

Namely, in the case of $3$-dimensional $0$-surgery, a \textit{horizontal} rotation of 180$^\circ$ around axis $\ell'$ transforms the two discs $S_1^0\times D^2$ of Fig.~\ref{2D_Decomp}({2\textsubscript{a}}) (the first instance  of $2$-dimensional $0$-surgery) to the two $3$-balls $S_1^0\times D^3$ of Fig.~\ref{3D_Decomp_F}(2\textsubscript{a}) (the first instance of $3$-dimensional $0$-surgery). After the collapsing of the centers of the two $3$-balls $S_1^0\times D^3$, the rotation transforms the decompactified cylinder $D^{1}\times (S_1^{1}\setminus \{\infty\})$ of Fig.~\ref{2D_Decomp}({2\textsubscript{a}}) (the last instance  of $2$-dimensional $0$-surgery) to the decompactified thickened sphere $D^{1}\times (S_2^{1}\setminus \{\infty\})$ of Fig.~\ref{3D_Decomp_F}(2\textsubscript{a}) (the last instance of $3$-dimensional $0$-surgery). Indeed, the rotation of line $\ell$ along $\ell'$ creates the green plane that cuts through ${\mathbb R}^3$ and separates the two resulting 3-balls $S_1^0\times D^3$. This plane is shown in green in the last instance of Fig.~\ref{3D_Decomp_F}(2\textsubscript{a}) and it is the decompactified view of the sphere $S_2^2$ in  ${\mathbb R}^3$. Note that it is thickened by the arcs connecting the two discs $S_1^0\times D^2$ which have also been rotated.

Similarly, in the case of $3$-dimensional $1$-surgery, a \textit{vertical} rotation of 180$^\circ$ around axis $\ell$ transforms the two discs $S_1^0\times D^2$ (the first instance  of $2$-dimensional $0$-surgery shown in Fig.~\ref{2D_Decomp}({2\textsubscript{a}})) to the solid torus $S_1^1\times D^2$ (the first instance of $3$-dimensional $1$-surgery), see Fig.~\ref{3D_Decomp_F}(3\textsubscript{a}). After the collapsing of the (red) core $S_1$ of $S_1^1\times D^2$, the rotation transforms the decompactified cylinder $D^{1}\times (S_1^{1}\setminus \{\infty\})$ of Fig.~\ref{2D_Decomp}({2\textsubscript{a}}) (the last instance  of $2$-dimensional $0$-surgery) to the decompactified solid torus $D^{2}\times (S_1^{1}\setminus \{\infty\})$ of Fig.~\ref{3D_Decomp_F}(3\textsubscript{a}) (the last instance of $3$-dimensional $1$-surgery). Indeed, each of the arcs $D^1$ connecting the two discs $S_1^0\times D^2$ generates through the rotation a $2$-dimensional disc $D^2$, and the set of all such discs are parametrized by the points of the line $\ell$ in ${\mathbb R}^3$.

In both cases, in Fig.~\ref{3D_Decomp_F}({1\textsubscript{a}}), $S^3$ is presented as the result of rotating the $2$-sphere $S^2 = {\mathbb R}^2 \cup \{\infty\}$. For $3$-dimensional $0$-surgery, $S^2$ is rotated about the circle $\ell' \cup \{\infty\}$  where $\ell'$ is a straight horizontal line in ${\mathbb R}^2$. The resulting decomposition of $S^3$ is $S^3 =(S_1^0\times D^3) \cup (D^1 \times S_2^2)$, a thickened sphere with two $3$-balls glued along the boundaries, which is visualized as $S^3\setminus \{\infty\} =(S_1^0\times D^3) \cup (D^1 \times (S_2^2\setminus \{\infty\}))$. For $3$-dimensional $1$-surgery, $S^2$ is rotated about the circle $\ell \cup \{\infty\}$  where $\ell$ is a straight vertical line in ${\mathbb R}^2$. The resulting decomposition of $S^3$ is $S^3 =(S_1^1\times D^2) \cup (D^2 \times S_2^1)$, two solid tori glued along their common boundary,  which is visualized as $S^3\setminus \{\infty\} =(S_1^1\times D^2) \cup (D^2 \times (S_2^1\setminus \{\infty\}))$. 
 
Analogously, starting with the second projection of Fig.~\ref{2D_Decomp}({1\textsubscript{b}}), the same rotations induce each type of $3$-dimensional surgery and their corresponding decompositions of $S^3$, see Fig.~\ref{3D_Decomp_F}({1\textsubscript{b}}). More precisely, a \textit{horizontal} rotation of the instances of Fig.~\ref{2D_Decomp}({2\textsubscript{b}}) by 180$^\circ$ around axis $\ell'$ induces the initial and final instances of $3$-dimensional $0$-surgery visualized in ${\mathbb R}^3$, see Fig.~\ref{3D_Decomp_F}({2\textsubscript{b}}). The $3$-sphere $S^3$ is now visualized as $S^3\setminus \{\infty\} =((S_1^0\setminus \{\infty\})\times D^3) \cup (D^1 \times S_2^2)$, a thickened sphere union two $3$-balls with the center of one of them removed (being the point at infinity).

Similarly, a  rotation of the instances of Fig.~\ref{2D_Decomp}({2\textsubscript{b}}) by 180$^\circ$ around the (green) \textit{circle} $\ell$ induces the initial and final instances of $3$-dimensional $1$-surgery visualized in ${\mathbb R}^3$, see Fig.~\ref{3D_Decomp_F}({3\textsubscript{b}}). Note that $\ell$ is now a circle and not a (vertical) line. The easiest part for visualizing this rotation is the rotation of the middle annulus of Fig.~\ref{3D_Decomp_F}({1\textsubscript{b}}) which gives rise to the solid torus $D^2 \times S_2^1$ in  Fig.~\ref{3D_Decomp_F}({3\textsubscript{b}}). The same rotation of the two remaining discs around $\ell$ can be visualized as follows: each radius of the inner disc lands from above the plane on the corresponding radius of the outer disc. At the same time, that radius of the outer disc lands on the corresponding radius of the inner disc from underneath the plane. So, the two corresponding radii together have created by rotation an annular ring around $\ell$. Note that the red center of the inner disc will land on all points at infinity, creating a half-circle from above and, at the same time, all points at infinity land on the center of the inner disc and create a half-circle from below. Glued together, the two half-circles create a (red) circle. Now, the set of all annular rings around $\ell$ and parametrized by $\ell$ make up the complement solid torus $S_1^1\setminus \{\infty\}\times D^2$ whose core is the aforementioned red circle. The $3$-sphere $S^3$ is visualized through this rotation as  $S^3\setminus \{\infty\} =((S_1^1\setminus \{\infty\})\times D^2) \cup (D^2 \times S_2^1)$, the decompactified union of two solid tori.

Finally, it is worth pinning down that the two types of visualizations presented above are related. Indeed, the $D^1 \times (S_2^2\setminus \{\infty\})$ shown in the rightmost instance of Fig.~\ref{3D_Decomp_F}({2\textsubscript{a}}) is the decompactified view of the $D^1 \times S_2^2$ shown in the rightmost instance of Fig.~\ref{3D_Decomp_F}({2\textsubscript{b}}). Likewise, the $D^2 \times (S_2^1\setminus \{\infty\})$ shown in the rightmost instance of Fig.~\ref{3D_Decomp_F}({3\textsubscript{a}}) is the decompactified view of the the solid torus $D^2 \times S_2^1$ shown in the rightmost instance of Fig.~\ref{3D_Decomp_F}({3\textsubscript{b}}). Further, the $(S_1^0\setminus \{\infty\})\times D^3$ and $(S_1^1\setminus \{\infty\})\times D^2$ shown in the leftmost instances of Fig.~\ref{3D_Decomp_F}({2\textsubscript{b}}) and ({3\textsubscript{b}}) are the decompactified views of $S_1^0\times D^3$ and $S_1^1\times D^2$ shown in the leftmost instances of Fig.~\ref{3D_Decomp_F}({2\textsubscript{a}}) and ({3\textsubscript{a}}) respectively.


\section{Conclusion} \label{Conclusion} 
In this paper we present how to detect the topological changes of $3$-dimensional surgery via the fundamental group and we provide a new way of visualizing its elementary steps. As this topological tool is used in both the classification of $3$-manifolds and in the description of natural phenomena, we hope that this study will help our understanding of the topological changes occurring in $3$-manifolds from a mathematical as well as a physical perspective.

\section*{Acknowledgments}
Antoniou's work was partially supported by the Papakyriakopoulos scholarship which was awarded by the Department of Mathematics of the National Technical University of Athens. Kauffman's work was supported by the Laboratory of Topology and Dynamics, Novosibirsk State University (contract no. 14.Y26.31.0025 with the Ministry of Education and Science of the Russian Federation).


\begin{thebibliography}{0}
\bibitem{Ad} Adams	C.:	The	Knot	Book,	An	Elementary	Introduction	to	the	Mathematical	Theory	of	Knots.	\textit{American	Mathematical	Society}	(2004).		

\bibitem{SS4}	Antoniou S., Kauffman L.H, Lambropoulou	S.: Topological surgery in cosmic phenomena (submitted for publication).		

\bibitem{BHsurg} Antoniou S., Kauffman L.H, Lambropoulou S.: Black holes and topological surgery. \textit{Preprint} (2018). Available from: https://arxiv.org/pdf/1808.00254.	
	
\bibitem{SS3}	Antoniou	S.:	Mathematical Modeling Through Topological Surgery and Applications. \textit{Springer Theses Book Series. Springer International Publishing}. DOI:  10.1007/978-3-319-97067-7 (2018).
				
\bibitem{SS1}	Antoniou	S.,	Lambropoulou	S.:	Extending	Topological	Surgery	to	Natural	Processes	and	Dynamical	Systems.	\textit{PLOS	ONE 12(9).}	DOI:	10.1371/journal.pone.0183993 (2017).										
							
\bibitem{SS2}	Antoniou	S.,	Lambropoulou	S.:	Topological	Surgery	in	Nature.	Book	`Algebraic	Modeling	of	Topological	and	Computational	Structures	and	Applications',	\textit{Springer	Proceedings	in	Mathematics	and	Statistics	Vol. 219}. DOI: 10.1007/978-3-319-68103-0 (2017).																						
\bibitem{StiGompf}	Gompf R., Stipsicz A.: 4-Manifolds and Kirby Calculus. \textit{Graduate Studies in Mathematics, Volume 20, American Mathematical Society} (1999).																										
\bibitem{DNA} Kauffman L.H, Lambropoulou S., Buck D.: DNA Topology. (book in preparation). Chapter on ‘Three Dimensional Topology for DNA’.							

\bibitem{KirScha}	 Kirby R.C., Scharlemann M.G.: Eight faces of the Poincar\'{e}  homology 3-sphere. \textit{Uspekhi Matematicheskikh Nauk, [N. S.]}. DOI: 10.1016/B978-0-12-158860-1.50015-0 (1979).																				
\bibitem{Kirby} Kirby R.: A Calculus for Framed Links in $S^3$, \textit{Inventiones math. 45, 35- 56} (1978).																							
\bibitem{Kosniowski} Kosniowski C.: A First Course in Algebraic Topology.  \textit{Cambridge: Cambridge University Press.} DOI: 10.1017/CBO9780511569296 (1980).	

\bibitem{SS5} Lambropoulou	S.,	Antoniou	S.:	Topological	Surgery,	Dynamics	and	Applications	to	Natural	Processes.	\textit{Journal	of	Knot	Theory	and	its	Ramifications 26(9)}.	DOI:	10.1142/S0218216517430027 (2016).																										
\bibitem{SS6} Lambropoulou	S.,	Samardzija	N.,	Diamantis	I.,	Antoniou	S.:	Topological	Surgery	and	Dynamics,	\textit{Mathematisches	Forschungsinstitut	Oberwolfach	Report	No.	26/2014,	Workshop:	Algebraic	Structures	in	Low-Dimensional	Topology}.	DOI:	10.4171/OWR/2014/26	(2014).	
							
\bibitem{Levin}	Levin  J.: Topology and the Cosmic Microwave Background. \textit{Physics Reports 365, 251–333} (2002).																										
\bibitem{LickTh} Lickorish W.B.R: A representation of orientable combinatorial 3-manifolds, \textit{Ann. of Math. (2) 76, pp. 531-540} (1962).																											
\bibitem{Luminet} Luminet, J.-P., Weeks J.R, Riazuelo A., Lehoucq R.,Uzan, J.-P.: Dodecahedral space topology as an explanation for weak wide-angle temperature correlations in the cosmic microwave background. \textit{Nature 425}, 593–595 (2003).																							\bibitem{MilAct}	Milnor	J.:	On the 3-dimensional Brieskorn manifolds M(p,q,r). \textit{Knot Groups and 3-Manifolds - Papers dedicated to the memory of R.H.Fox. Annals of Mathematics Studies 84. Princeton University Press, Princeton, NJ} (1975).				
\bibitem{Milsur} Milnor J.:	A procedure for killing the homotopy groups of differentiable manifolds. \textit{Symposia in Pure Math., Amer. Math. Soc. 3, 39-55} (1961).																												
\bibitem{Munkres} Munkres J.R: Topology, Second Edition. \textit{Prentice Hall, Incorporated} (2000).																											
\bibitem{PS}	Prasolov	V.V.,	Sossinsky,	A.B.:	Knots,	links,	braids	and	3-manifolds.	\textit{AMS	Translations	of	Mathematical	Monographs	154}	(1997).
									
\bibitem{Ra}	Ranicki	A.:	Algebraic	and	Geometric	Surgery.	\textit{Oxford	Mathematical	Monographs,	Clarendon	Press}	(2002).															
\bibitem{Rolfsen}  Rolfsen	D.:	Knots	and	links.	\textit{Publish	or	Perish	Inc.	AMS	Chelsea	Publishing}	(2003).															
\bibitem{N7}	Samardzija	N.,	Greller	L.:	Explosive	route	to	chaos	through	a	fractal	torus	in	a	generalized	Lotka-Volterra	model.	\textit{Bull	Math	Biol	50,	No.	5:00	465--491}.	DOI:	10.1007/BF02458847	(1988).
\bibitem{Sti}	Stillwell J.:	Classical Topology and Combinatorial Group Theory. \textit{Graduate Texts in Mathematics, Springer-Verlag New York} (1993).																									
\bibitem{SeTh}	Threlfall H., Seifert W. : A textbook of topology. \textit{Academic Press} (1980).																										
\bibitem{Wal}	Wallace A.H.:	Modifications and cobounding manifolds. \textit{Canad. J.  Math. 12, 503-528} (1960).																									
\bibitem{Weeks}	Weeks  J.R.: The Shape of Space. \textit{CRC Press} (2001).																										
\bibitem{MagnusKarassSolitar} Wilhelm M., Karrass A., Solitar D.: Combinatorial Group Theory, \textit{New York: Dover Publications} (2004).																											
\end{thebibliography}
\end{document}